\documentclass [12pt] {amsart}
\usepackage{amsthm,amsmath,amssymb,amscd,graphics,enumerate}
\usepackage{latexsym}
\usepackage{verbatim}
\usepackage{enumerate}
\usepackage{amsthm}
\usepackage{amsmath}
\usepackage{amscd}

\newcounter{theorem}[section]
\numberwithin{equation}{section}

\newtheorem{Main}[theorem]{Main Theorem}
\newtheorem{Thm}[theorem]{Theorem}

\newtheorem{lemma}[subsection]{Lemma}

\newtheorem{proposition}[subsection]{Proposition}
\theoremstyle{definition}

\theoremstyle{remark}

\theoremstyle{remark}

\numberwithin{equation}{section}
{\theoremstyle{remark}
 
 }

{\theoremstyle{definition}
 }

\newcommand{\RR}{{\mathbb{R}}}

\newcommand{\PP}{{\mathbb{P}}}

\newcommand{\mO}{\mathcal{O}}
\newcommand{\mX}{\mathcal{X}}
\newcommand{\mY}{\mathcal{Y}}
\newcommand{\mW}{\mathcal{W}}

\newcommand{\mC}{\mathcal{C}}

\newcommand{\mD}{\mathcal{D}}
\newcommand{\mE}{\mathcal{E}}
\newcommand{\mF}{\mathcal{F}}

\newcommand{\mH}{\mathcal{H}}

\newcommand{\mV}{\mathcal{V}}

\newcommand{\X}{\mathcal{X}}

\newcommand{\D}{\mathcal{D}}

\newcommand{\sO}{\mathcal{O}}

\def\<{\left\langle}
\def\>{\right\rangle}

\begin{document}

\title[Stacky Bogomolov-Miyaoka-Yau inequality]{On the Bogomolov-Miyaoka-Yau inequality\\ for stacky surfaces}
\author{Jiun-Cheng Chen}
\address{Department of Mathematics\\ Third General Building, National Tsing Hua University\\ No. 101 Sec 2 Kuang Fu Road \\Hsinchu, Taiwan 30043, Taiwan}
\email{jcchen@math.nthu.edu, jcchenster@gmail.com}

\author{Hsian-Hua Tseng}
\address{Department of Mathematics\\ Ohio State University\\ 100 Math Tower, 231 West 18th Ave.\\Columbus\\ OH 43210\\ USA}
\email{hhtseng@math.ohio-state.edu}

\date{\today}

\begin{abstract}
We discuss a generalization of the Bogomolov-Miyaoka-Yau inequality to Deligne-Mumford surfaces of general type.
\end{abstract}
\maketitle
\section{Introduction}
We work over $\mathbb{C}$. 

For a smooth complex projective surface $S$ of general type, the Bogomolov-Miyaoka-Yau inequality for $S$ reads (see \cite{mi77})
\begin{equation}\label{BMY_sm_surface}
3c_2(T_S)\geq c_1(T_S)^2.
\end{equation}
Together with Noether's inequality, this puts constraints on the topology of surfaces of general types. Generalizations of (\ref{BMY_sm_surface}) to singular surfaces and surface pairs have been found, see for example \cite{m}, \cite{l1, l2}. 

In this paper we discuss a generalization of (\ref{BMY_sm_surface}) to Deligne-Mumford stacks. Let $\mX$ be a smooth proper Deligne-Mumford $\mathbb{C}$-stack of dimension $2$. Let $\pi:\mX \to X$ be the natural map to the coarse moduli space. We assume that $X$ is a projective variety. Since $\mX$ is assumed to be smooth, it has a tangent bundle $T_\mX$. A good theory of Chern classes is available for Deligne-Mumford stacks, see for example \cite{v}, \cite{k}. 

\begin{Main}
Let $\mX$ be as above. Assume that the canonical bundle $K_\mX:=\wedge^2T_\mX^\vee$ is numerically effective, then
\begin{equation}\label{BMY_stack}
3 c_2(T_\mX)\geq c_1(T_\mX)^2.
\end{equation}
\end{Main}

Certainly (\ref{BMY_stack}) takes the same shape as (\ref{BMY_sm_surface}). A proof of (\ref{BMY_stack}), along the lines of Miyaoka's original proof of (\ref{BMY_sm_surface}) in \cite{mi77}, is given in Section \ref{pf_BMY}. Section \ref{sec:examples} contains examples of (\ref{BMY_stack}). In Section \ref{codim1} we consider (\ref{BMY_stack}) for a class of stacks $\mX$ with stack structures in codimension $1$, recovering \cite[Corollary 0.2]{l1}. In Section \ref{codim2} we consider (\ref{BMY_stack}) for Gorenstein stacks $\mX$ with isolated stack points, recovering \cite[Corollary 1.3]{m}. 

Generalizations of the Bogomolov-Miyaoka-Yau inequality to varieties with quotient singularities (i.e. orbifolds) certainly have been studied before by many authors using various approaches. References to these can be found in e.g. \cite{l1, l2}. In this paper we work in the context of Deligne-Mumford stacks. This viewpoint has the advantage that (\ref{BMY_stack}) can be proven by following Miyaoka's original arguments in \cite{mi77}. 
Also, as discussed in Section \ref{sec:examples}, (\ref{BMY_stack}) specializes to some generalizations of the original (\ref{BMY_sm_surface}) by straightforward and elementary means.  

 \subsection*{Acknowledgment} 
J. -C. C. is a Golden-Jade Fellow of Kenda Foundation, Taiwan. He is supported in part by National Science Council and National Center for Theoretical Sciences, Taiwan. H.-H. T. is supported in part by NSF grant DMS-1047777 and Simons Foundation Collaboration Grant.

\section{Proof of (\ref{BMY_stack})}\label{pf_BMY}
In this Section we give a proof of (\ref{BMY_stack}). Our proof is adapted from Miyaoka's original proof in \cite{mi77}.

Let $\mX$ be a smooth proper Deligne-Mumford stack of dimension $2$. If $\X$ has non-trivial stack structures at generic points, then $\X$ is an \'etale gerbe over a stack with trivial generic stack structure, see for example \cite[Proposition 4.6]{bn}. More precisely, there is a finite group $G$, a stack $\X'$ with trivial generic stabilizers, and a morphism $f: \X\to \X'$ realizing $\X$ as a $G$-gerbe over $\X'$. Since $T_\X=f^*T_{\X'}$, we see that (\ref{BMY_stack}) for $\X$ is equivalent to (\ref{BMY_stack}) for $\X'$. Therefore it suffices to consider only those $\X$ with stack structures in codimension $\geq 1$.  For the rest of this section we assume this.

Let $\mF$ be a locally free sheaf of rank $2$ on $\X$. Let $\mV:=\mathbb{P}(\mF)$ be the projectivization, with natural projection $p: \mV\to \X$. Let $\mH$ be the divisor associated to the tautological sheaf on $\mV$. 

\begin{lemma}\label{fundamentallemma}
Assume that $\mW \subset \mV$ is linearly equivalent to $\mH- p^* \mD$, where $\mD \subset \mX$ is a divisor on $\mX$. Then we have
\begin{equation}
\mD \cdot det \mF \leq c_2(\mF) + \mD ^2.
\end{equation}
\end{lemma}

\begin{proof}
We closely follow Miyaoka's original proof \cite{mi77}. Let $i: \mW \subset \mV$ be the inclusion morphism. Note that
the composition $p \circ i: \mW \to \mX$ is birational by our assumption on the linear equivalence class of $\mW$.
Since resolutions can be chosen such that they are compatible with \'etale base change, 
there is a sequence of blow-ups 
\begin{equation}
\mu: \mV_s  \xrightarrow{\mu_s}  \mV_{s-1} \rightarrow \cdots \rightarrow \mV_1  \xrightarrow{\mu_1}  \mV_0=\mV
\end{equation}
such that the proper transform $ \mW'$ of $\mW$ is a smooth Deligne-Mumford stack in $\mV_s$. Let $i': \mW'\subset \mV_s$ and $\rho: \mW'\to \X$ be the natural maps.

Let $\mE_1$, $\mE_2$, $\cdots$, $\mE_s$ be the exceptional divisors on $\mV_s$. The divisor $\mW'$ is linearly equivalent to 
$\mu ^* (\mH- p^*\mD) -\sum a_i \mE_i$. It can be seen\footnote{The argument is similar to that of \cite[Lemma 7]{mi77} and is omitted.} that the canonical bundle $K_{\mW'}$ satisfies $K_{\mW'}= \rho^*K_{\mX}+ \sum \mC_i$ where $C_i$ is a curve and $\rho( \mC_i)=point$.
By the Hodge index theorem (for a stacky version see \cite[Theorem 3.1.3]{lie}), it follows that
$(K_{\mW'}- \rho ^*K_{\mX} + \sum c_i i'^*\mE_i)^2 \leq 0$ for any $c_i \in \RR$.

Write $K_{\mV_{s}}= \mu^*(-2\mH+ p^* K_{\mX}+ p^*(det \mF))+ \sum b_i \mE_i$. 
The adjunction formula implies that  
$$K_{\mW'}= i'^*[\mu^*(-\mH)+ (p \circ \mu)^*(K_{\mX}+det \mF -\mD) + \sum (b_i-a_i) \mE_i].$$
Thus $i'^*\mu^*(-\mH+  p^*(det \mF -\mD))^2 \leq 0$.
Set $k:= i'^* \mu^*(-\mH+  p^*(det \mF -\mD))^2$. 
We can also compute this self-intersection number $k$ in another way:
\begin{equation*}
\begin{split}
k&= \mu^*(- \mH+ p^*(det \mF)-p^* \mD)^2 ( \mu^* \mH- ( p \circ \mu )^* \mD -\sum a_i \mE_i)\\
&= \mu^*(- \mH+ p^*(det \mF)-p^* \mD))^2 ( \mu^* \mH- ( p \circ \mu )^* \mD) \quad \quad \quad (\text{since $\mE_i$ is exceptional})\\
&=\mH^3-\mH^2 \cdot p^*(\mD -2 det \mF)+ \mH \cdot (p^*(det \mF)^2-(p^* \mD)^2).
\end{split}
\end{equation*}
Using the standard relations for the intersection numbers on the
 projectivization of a rank $2$ vector bundle, we calculate that
\begin{equation*}
k= c_1^2(det(\mF))-c_2(\mF) -(det \mF)^2+ det \mF \cdot \mD -\mD^2= -c_2(\mF) +det \mF \cdot \mD -\mD^2.
\end{equation*}
The result follows.
\end{proof}

Let $\mO_{\mX}( \mD)$ be a subsheaf of $\Omega^1_{\mX}$. One key observation used in Miyaoka's original proof is that the Iitaka dimension 
of $\mO_{\mX}( \mD)$ is at most $1$.
 
\begin{Thm}\label{iitakadim}
If $\mO_{\mX}( \mD)$ is a subsheaf of $\Omega^1_{\mX}$ of a projective Deligne-Mumford stack $\mX$, then 
$h^0(\mX, \mO_{\mX}(n \mD)) \leq c n$ for some positive constant $c$ and $n>>0$. 
\end{Thm}

The proof of Theorem \ref{iitakadim} is very similar to that of \cite[Theorem 2'']{mi77}.
Two main ingredients are needed in the proof of \cite[Theorem 2'']{mi77}: (1) the Riemann-Roch formula, and (2) a lemma due to de Franchis.  The de Franchis lemma states that any global holomorphic differential form  on a K\"{a}hler manifold or a surface is $d$-closed \cite[Lemma 9]{mi77}. This lemma follws essentially from  Stoke's theorem. The  argument still works for  smooth K\"{a}hler Deligne-Mumford stacks (K\"{a}hler orbifolds) or smooth surface Deligne-Mumford stacks. The Riemann-Roch for stacks is proved in \cite{toen_rr}.

One can prove the following result using Theorem \ref{iitakadim}. 
\begin{proposition}\label{keyprop}
Let $ \mF \subset \Omega ^1_{\mX}$ be a locally free sheaf of rank 2 and assume that $det(\mF) ^{\otimes n}$ is generated by global sections for some $n>0$. If $\mF \otimes \mO_{\mX}(- \mD)$ has a non-trivial section, then 
\begin{equation*}
\mD \cdot det(\mF) \leq max \{ c_2(\mF),0\}.
\end{equation*}
\end{proposition}

\begin{proof}
Consider $p: \mV=\PP(\mF) \to \mX$. The canonical isomorphism gives us $$H^0(\mX, \mF \otimes \mO_{\mX}(- \mD))=H^0(\mV, \mO_{\mV}(\mH -p^* \mD)).$$
If $\mF \otimes \mO_{\mX}(- \mD)$ has a non-trivial section, then $|\mH -p^* \mD|$ is non-empty. Pick $\mW \in |\mH -p^* \mD|$. Decompose
$\mW$ as $\mW= \mW_0+ p^* \mD'$ where $\mW_0$ is effective and irreducible which is linearly equivalent to $\mH- p^*(\mD+\mD')$ and $\mD'$ is effective. 
Note that $(det \mF)^{\otimes n}$ is generated by global sections, so the intersection number $\mD' \cdot det (\mF) \geq 0$.  It follows that 
$\mD \cdot det (\mF) \leq (\mD + \mD') \cdot det (\mF)$ and it suffices to prove $(\mD+\mD') \cdot det (\mF) \leq max \{c_2(\mF),0\}$. Set $\mD''= \mD+\mD'$ to simplify notation.  By Lemma \ref{fundamentallemma}, $\mD'' \cdot det(\mF) \leq c_2(\mF) + \mD'' \cdot \mD''$. Observe that $\mO_{\mX}(\mD'')$ is a subsheaf of $\Omega^1_{\mX}$.  Indeed, the effectiveness of $\mW_{0}$ ensures the existence  of  a non-trivial section
 of $\mF \otimes \mO_{\mX}(- \mD'')$, i.e. an injection $ \mO_{\mX} \hookrightarrow  \mF \otimes \mO_{\mX}(- \mD'')$.  Twisting by $\mO_{\mX}(- \mD'')$, embeds $\mO_{\mX}(- \mD'')$ into $\mF \subset \Omega^{1}_{\mX}$.  By Theorem \ref{iitakadim},   $\mD''$ has Iitaka dimension at most $1$. 
It follows that $\mD'' \cdot det(\mF) \leq 0$ or $\mD'' \cdot \mD'' \leq 0$. \footnote{Arguing as in  \cite[Lemma 10]{mi77}.} 
This completes the proof.
\end{proof}

Assuming $c_2(\mF)$ is positive for the time being, we can obtain an upper bound on $c_2$ provided the sheaf $\mF \otimes \mO_{\mX}(- \mD)$ has no sections. This can then be used to derive a contradiction. To be more precise, one needs a  modified version of Proposition \ref{keyprop}, in which the condition on the sheaf $\mF \otimes \mO_{\mX}(- \mD)$ having a section is replaced by the condition that some symmetric power $S^m \mF \otimes \mO_{\mX}(- \mD)$ having a section.

\begin{Thm}\label{boundonc2}
Let $ \mF \subset \Omega ^1_{\mX}$ be a locally free sheaf of rank 2 and assume that $det(\mF) ^{\otimes n}$ is generated by global sections for some $n>0$. If $S^m \mF \otimes \mO_{\mX}(- \mD)$ has a non-trivial section, then 
\begin{equation*}
\mD \cdot det(\mF) \leq max \{ mc_2(\mF),0\}.
\end{equation*}
\end{Thm}

The proof of Theorem \ref{boundonc2} follows from Proposition \ref{keyprop} and  the following easy lemma (which is analogous to \cite[Lemma 11]{mi77}).

\begin{lemma}\label{basechange}
Let $p: \mV= \PP(\mF) \to \mX$ be the projective bundle of a locally free sheaf of rank $2$. Let $\mW \in | m\mH- p^* \mD|$.
Then there is a surjective morphism $\beta: \mX' \to \mX$ such that $\beta ^*\mW$ is decomposed to $\mW_1+ \cdots \mW_m$ where
$\mW_i$ is an effective divisor linear equivalent to $\mH' -p^* \mD_i$. 
\end{lemma}

\begin{proof}[Proof of Theorem \ref{boundonc2}]
The following argument is taken from \cite[Theorem 3]{mi77}. Let $f$ be a global section of $S^m \mF \otimes \mO_{\mX}(- \mD)$. 
Lemma \ref{basechange} imples that after a suitable cover $\beta:\mY \to \mX$, we can decompose $\beta ^* f$ can be written as  $f_1 f_2 \cdots f_m \in H^0(\mY, S^m \beta^{*} \mF \otimes \mO_{\mY}( -\beta ^* \mD))$, where $f_i \in H^0(\mY, \beta^{*} \mF \otimes \mO_{\mY}( -\beta ^*\mD_{i}))$ and $(det \beta^* \mF)^{\otimes m} \cong (\beta^* det ( \mF)^{\otimes m})$ is generated by global sections. From Proposition \ref{keyprop}, it follows that 
$\beta^* \mD_i \cdot (det( \beta ^*\mF)) \leq max \{ c_2(\beta^* \mF),0\}$. Summing over all $i$'s, we have $\beta ^*\mD \cdot  det(\beta^* \mF) \leq max \{m c_2(\beta^*\mF),0\}.$ Let $d$ be the mapping degree of $\beta$. Clearly, $\beta^* \mD \cdot det( \beta ^*\mF)=d D \cdot det(\mF)$ and $c_2( \beta^* \mF)=d \beta ^*c_2(\mF)$.
 \end{proof}

We now come to (\ref{BMY_stack}).
\begin{Thm}\label{mainthm}
Let $\mX$ be a non-singular Deligne-Mumfors stack with  the projective coarse space $X$ of general type and $c_1(\mX)$ nef. Then 
$c_1^2(\mX) \leq 3 c_2(\mX)$ holds.  
\end{Thm}

\begin{proof}

As in \cite{mi77}, we consider two cases: (1) $c_1^2(\mX) \leq 2 c_2(\mX)$ and (2) $c_1^2(\mX) >2 c_2(\mX)$. 
The first case is obvious. For the second case, set $\alpha := \frac{c_2(\mX)}{c_1^2(\mX)}$. Note that $\alpha <1/2$. Pick $\delta >0$ sufficiently small and rational. By Theorem \ref{boundonc2} applied to $\mD=m(\alpha+\delta)K_\X$, $\mF=\Omega_\X^1$, 
we can find a positive integer $m$ such that $m (\alpha + \delta) \in \mathbb{Z}$, and
$$ h^0( \mX, S^m \Omega_{\mX}^1 \otimes \mO_{\mX}(-m(\alpha + \delta) K_{\mX}))=0.$$
By Serre duality for smooth projective Deligne-Mumford stacks \cite[Theorem 2.22]{ni08}, we have 
$$h^2( \mX, S^m \Omega_{\mX}^1 \otimes \mO_{\mX}(-m(\alpha + \delta) K_{\mX}))=h^0( \mX, S^m \Omega_{\mX}^1 \otimes \mO_{\mX}(-m(1-\alpha - \delta) K_{\mX})\otimes K_{\mX}).$$
As $\alpha < 1/2$ and $\delta$ is small, we have $1- \alpha - \delta > \alpha$. We apply Theorem \ref{boundonc2} to $\mD=m(2-\alpha-\delta)K_\X$, $\mF=\Omega_\X^1$, to get $$h^2( \mX, S^m \Omega_{\mX}^1 \otimes \mO_{\mX}(-m(\alpha + \delta) K_{\mX}))=0.$$ 
Hence $$\chi ( \mX, S^m \Omega_{\mX} \otimes \mO (-m (\alpha+ \delta) K_{\mX}))= -h^1( \mX, S^m \Omega_{\mX}^1 \otimes \mO_{\mX}(-m(\alpha + \delta) K_{\mX})) \leq 0.$$

Note that to compute the cohomology groups of a (subsheaf of) symmetric power of a vector bundle, one can work on the the projectivized vector bundle and computing the  cohomology groups of relevant line bundles. Thus  $$0\geq \chi ( \mX, S^m \Omega_{\mX} \otimes \mO (-m (\alpha+ \delta) K_{\mX}))= \chi ( \mV, \mO_{\mV} (-m (\mH-(\alpha+ \delta) \pi^* K_{\mX})).$$
By Riemann-Roch for stacks \cite{toen_rr}, we have $\chi ( \mV, \mO_{\mV} (-m (\mH-(\alpha+ \delta) \pi^* K_{\mX}))$ grows like $\frac{1}{6} (\mH-(\alpha+ \delta) \pi^* K_{\mX})^3 m^3$ as $m\to \infty$.
It implies that $(\mH-(\alpha+ \delta) \pi^* K_{\mX})^3 \leq0$. Taking $\delta $ to $0$, we obtain 
\begin{equation*}
\begin{split}
0 \geq (\mH-\alpha \pi^* K_{\mX})^3
=&c_1^2(\mX)-c_2(\mX) -3 \alpha c_1^2(\mX)+ 3 \alpha ^2 c_1^2(\mX) \\
=&(1-\alpha -3 \alpha + 3 \alpha ^2) c_1 ^2(\mX)\\
=&(1-\alpha)(1-3\alpha)c_1^2(\mX).
\end{split}
\end{equation*}
Since $\alpha <1/2$ and $c_1 ^2(\mX)$ is non-negative, we get $1- 3 \alpha \leq 0$ as desired.
\end{proof}

\section{Examples of (\ref{BMY_stack})}\label{sec:examples}

\subsection{General discussion}
According to the main result of \cite{gs}, the map $\X\to X$ from the stack $\X$ to its coarse moduli space $X$ (which we assume to be a variety with quotient singularities) can be factored as $$\X\to \X_1\to \X_2\to X,$$
where 
\begin{enumerate}
\item
$\X_1$ has trivial generic stabilizers;
\item
$\X_2$ is the canonical stack associated to the variety $X$ (see e.g. \cite[Definition 4.4]{fmn}) and has stack structures in codimension at least $2$;
\item
$\X\to \X_1$ is a {\em gerbe};
\item 
$\X_1\to \X_2$ is a composition of {\em root constructions} along divisors, thus introducing codimension-$1$ stack structures to $\X_2$.
\end{enumerate}

Since $\X\to \X_1$ is a gerbe, the tangent bundle of $\X_1$ pulls back to the tangent bundle of $\X$. So the inequality (\ref{BMY_stack}) for $\X$ is equivalent to (\ref{BMY_stack}) for $\X_1$. Therefore when considering examples, we may restrict our attention to $\X$ whose stack structures are in codimension at least $1$. In the rest of this section we present two examples of (\ref{BMY_stack}): the example in Section \ref{codim1} is obtained by root constructions, and the examples in Section \ref{codim2} are canonical stacks associated to quotient varieties. In these examples we show that (\ref{BMY_stack}) coincides with previous results.

\subsection{Codimension $1$ stack structure}\label{codim1}
We consider (\ref{BMY_stack}) for an example of stack $\X$ with stack structures in codimension $1$. 

Let $X$ be a smooth complex projective surface and $D$ a simple normal crossing $\mathbb{Q}$-divisor of the form $D=\sum_i (1-1/r_i)D_i$ with $r_i\geq 2$ integers. Let $\X$ be the natural stack cover of the pair $(X, D)$. By construction the coarse moduli space of $\X$ is $X$. The natural map $\pi: \X\to X$ is an isomorphism outside $\pi^{-1}(\text{Supp}\,D)$, which is where $\X$ has non-trivial stack structures. The stack $\X$ can be constructed from $X$ by applying root constructions along components of $D$. Furthermore we have the following formula for the canonical bundle: 
\begin{equation}\label{K-formula}
K_\X=\pi^*(K_X+D).
\end{equation}

We now examine (\ref{BMY_stack}) for this $\X$. By (\ref{K-formula}), $$c_1(T_\X)^2=c_1(K_\X)^2=(K_X+D)^2.$$ By Gauss-Bonnet theorem for Deligne-Mumford stacks \cite[Corollaire 3.44]{toen} we have $$c_2(T_\X)=\chi(\X),$$ the Euler characteristic of $\X$ as defined in \cite[Definition 3.43]{toen} (note that the notation $\chi^{orb}$ is used in \cite{toen}). Put $$\D_i:=\pi^{-1}(D_i), \quad \D_i^\circ:=\D_i\setminus (\cup_{j\neq i}(\D_i\cap \D_j)).$$
Then we have 
\begin{equation*}
\chi(\X\setminus \pi^{-1}(\text{Supp}\,D))=\chi(\X)-\sum_i \chi(\D_i^\circ)-\sum_{p\in \D_i\cap \D_j}\chi(p).
\end{equation*}
Similarly, put $D_i^\circ=D_i\setminus (\cup_{j\neq i} (D_i\cap D_j))$, we have 
\begin{equation*}
\chi(X\setminus \text{Supp}\,D)=\chi(X)-\sum_i\chi(D_i^\circ)-\sum_{\bar{p}\in D_i\cap D_j} \chi(\bar{p}).
\end{equation*}
Since $\X\setminus\pi^{-1}(\text{Supp}\,D)\simeq X\setminus \text{Supp}\,D$, we have $\chi(\X\setminus \pi^{-1}(\text{Supp}\,D))=\chi(X\setminus \text{Supp}\,D)$. Equivalently, 
\begin{equation*}
\chi(\X)=\chi(X)-\sum_i\chi(D_i^\circ)-\sum_{\bar{p}\in D_i\cap D_j} \chi(\bar{p})+ \sum_i \chi(\D_i^\circ)+\sum_{p\in \D_i\cap \D_j}\chi(p).
\end{equation*}
Since the map $\D_i^\circ\to D_i^\circ$ is of degree $1/r_i$ and the map $\D_i\cap\D_j\to D_i\cap D_j$ is of degree $1/r_ir_j$, we have $$\chi(\D_i)=\frac{1}{r_i}\chi(D_i), \quad \chi(\D_i\cap \D_j)=\frac{1}{r_ir_j}\chi(D_i\cap D_j).$$
This implies that 
\begin{equation}\label{eulerstackX}
\chi(\X)=\chi(X)-\sum_i(1-1/r_i)\chi(D_i^\circ)+\sum_{\bar{p}\in D_i\cap D_j}(1/r_ir_j-1).
\end{equation}
By \cite[Theorem 8.7]{l2}, for $\bar{p}\in D_i\cap D_j$ the local orbifold Euler number of the pair $(X,D)$ at $\bar{p}$ is given by $e_{orb}(\bar{p}; X, D)=1/r_ir_j$. Together with (\ref{eulerstackX}) this implies that $\chi(\X)$ coincides with the orbifold Euler number $e_{orb}(X,D)$ of the pair $(X,D)$, as defined in \cite{l2}. Thus if $K_\X$ is numerically effective, then (\ref{BMY_stack}) is equivalent to \cite[Theorem 0.1]{l2} applied to the pair $(X,D)$.

\subsection{Codimension $2$ stack structure}\label{codim2}
Let $\X$ be a smooth proper Deligne-Mumford $\mathbb{C}$-stack of dimension $2$ with isolated stack structures. 
Let $\pi: \X\to X$ be the natural map to the coarse moduli space $X$. Let $p_1, p_2,...,p_k\in \X$ be the stacky points.  Suppose that $\X$ is Gorenstein, i. e. 
 each stacky point $p_i$ has a neighborhood $p_i\in U_i\subset \X$ of the form $U_i\simeq [\mathbb{C}^2/G_i]$ with $G_i\subset SU(2)$ a finite subgroup, identifying $p_i$ with $[0/G_i]\in [\mathbb{C}^2/G_i]$.  It is a  standard  fact that the coarse moduli space $X$ is a projective surface with canonical singularities.

Suppose further that $K_\X$ is numerically effective. We consider (\ref{BMY_stack}) for such $\X$. 

By assumption we have $K_\X=\pi^*K_X$. Thus $$c_1(T_\X)^2=c_1(K_\X)^2=c_1(K_X)^2.$$ We now consider the term $c_2(T_\X)$. The first step is to consider $\chi(\sO_\X)$ by using Riemann-Roch theorem for stacks \cite{toen, toen_rr}. We follow \cite[Appendix A]{t} for the presentation of the Riemann-Roch theorem. We have $$\chi(\sO_\X)=\int_{I\X}\widetilde{ch}(\sO_\X)\widetilde{Td}(T_\X).$$ Here $I\X$ is the inertia stack of $\X$. By our assumption on $\X$, we have the following description of $I\X$:
$$I\X=\X\cup \bigcup_{i=1}^k (Ip_i \setminus p_i).$$ Here the term $Ip_i\setminus p_i$ is the inertia stack of $p_i\simeq BG_i$ with the main component removed, namely $$Ip_i\setminus p_i\simeq \bigcup_{(g)\neq (1): \text{conjugacy class of } G_i} BC_{G_i}(g).$$
By the definition of the Chern character $\widetilde{ch}$, we have $\widetilde{ch}(\sO_\X)=1$ on every component of $I\X$. Hence 
\begin{equation}\label{RR-sum}
\chi(\sO_\X)=\int_{I\X}\widetilde{Td}(T_\X)=\int_\X\widetilde{Td}(T_\X)|_\X+\sum_{i=1}^k\int_{Ip_i\setminus p_i} \widetilde{Td}(T_\X)|_{Ip_i\setminus p_i}.
\end{equation}
Note that $\widetilde{Td}(T_\X)|_\X=Td(T_\X)$, and we only need its degree $2$ component. Hence 
\begin{equation}\label{RR-main-term}
\int_\X\widetilde{Td}(T_\X)|_\X=\frac{1}{12}\int_\X (c_2(T_\X)+c_1(T_\X)^2).
\end{equation}
The contribution coming from $Ip_i\setminus p_i$ can be also evaluated.
\begin{lemma}\label{RR-computation}
Let $E_i$ be the exceptional divisor of the minimal resolution of $\mathbb{C}^2/G_i$. Then $$\int_{Ip_i\setminus p_i} \widetilde{Td}(T_\X)|_{Ip_i\setminus p_i}=\frac{1}{12}(\chi(E_i)-\frac{1}{|G_i|}).$$
\end{lemma}
An elementary proof of this Lemma is given in the Appendix. 

Next, we reinterpret the term $\chi(\sO_\X)$. By definition, $\chi(\sO_\X):=\sum_{l\geq 0}(-1)^l \text{dim}\, H^l(\X, \sO_\X)$. Since $\pi_*\sO_\X=\sO_X$ (see e.g. \cite[Theorem 2.2.1]{av}), we have $H^l(\X, \sO_\X)=H^l(X, \sO_X)$ and 
\begin{equation}\label{euler-equality}
\chi(\sO_\X)=\chi(\sO_X).
\end{equation}

Combining (\ref{RR-sum}), (\ref{RR-main-term}), (\ref{euler-equality}), and Lemma \ref{RR-computation}, we obtain the following expression of $c_2(T_\X)$:
\begin{equation}\label{2nd_chern}
\int_\X c_2(T_\X)=12\chi(\sO_X)-\int_\X c_1(T_\X)^2-\sum_{i=1}^k(\chi(E_i)-1/|G_i|).
\end{equation}
Using this, we see that in the present situation, (\ref{BMY_stack}) is equivalent to 
\begin{equation}\label{BMY_stack_codim2}
12\chi(\sO_X)\geq \frac{4}{3}c_1(K_X)^2+\sum_{i=1}^k(\chi(E_i)-\frac{1}{|G_i|}).
\end{equation}
On the other hand, it is clear that (\ref{BMY_stack_codim2}) is a special case of \cite[Corollary 1.3]{m}.

\appendix
\section{Proof of Lemma \ref{RR-computation}}
In this Appendix we prove Lemma \ref{RR-computation}. By our assumption on $\X$, for $g\in G_i$, the $g$-action on the tangent space $T_{p_i}\X$ has two eigenvalues $\xi_g$ and $\xi_g^{-1}$, where $\xi_g$ is a certain root of unity. By the definition of $\widetilde{Td}(T_\X)$ we have 
\begin{equation}\label{twisted_sector_contribution}
\int_{Ip_i\setminus p_i} \widetilde{Td}(T_\X)|_{Ip_i\setminus p_i}=\sum_{(g)\neq (1): \text{conjugacy class of } G_i}\frac{1}{|C_{G_i}(g)|}\frac{1}{2-\xi_g-\xi_g^{-1}}.
\end{equation}
We now evaluate (\ref{twisted_sector_contribution}) using the $ADE$ classification of $\mathbb{C}^2/G_i$. 

\subsection{Type A}
If $\mathbb{C}^2/G_i$ is of type $A_{n-1}$, then $G_i\simeq \mathbb{Z}_n$ and the action on $\mathbb{C}^2$ is given as follows. If we identify $\mathbb{Z}_n$ with the group of $n$-th roots of $1$, then an element $\xi\in \mathbb{Z}_n$ acts on $\mathbb{C}^2$ via the matrix 
$$\left(\begin{array}{cc}
\xi& 0 \\
 0& \xi^{-1}
 \end{array}\right).$$
It follows that (\ref{twisted_sector_contribution}) is given by 
\begin{equation}\label{cont_A}
\frac{1}{n}\sum_{l=1}^{n-1}\frac{1}{2-\exp(2\pi\sqrt{-1}l/n)-\exp(2\pi\sqrt{-1}l/n)^{-1}}.
\end{equation}
By \cite[Lemma 3.3.2.1]{lie}, (\ref{cont_A}) is equal to $$\frac{n^2-1}{12n}=\frac{1}{12}(n-1/n).$$ Since the exceptional divisor of the minimal resolution of $\mathbb{C}^2/\mathbb{Z}_n$ is a chain of $(n-1)$ copies of $\mathbb{CP}^1$, its Euler characteristic is $n$. This proves the Lemma in type A case.

\subsection{Type D}
If $\mathbb{C}^2/G_i$ is of type $D_{n+2}$ (here $n\geq 2$), then $G_i$ is isomorphic to the binary dihedral group $Dic_n$. The group $Dic_n$ is of order $4n$ and may be presented as follows: $$Dic_n=\<a,x|a^{2n}=1, x^2=a^n, x^{-1}ax=a^{-1}\>.$$ The action of $Dic_n$ on $\mathbb{C}^2$ is given as follows: 
\begin{equation}\label{action_typeD}
a\mapsto  \left(\begin{array}{cc} \exp(\pi\sqrt{-1}/n)& 0\\ 0& \exp(-\pi\sqrt{-1}/n)\end{array}\right), \quad x\mapsto \left(\begin{array}{cc} 0& 1\\ -1& 0\end{array}\right).
\end{equation}
An elementary calculation shows that the conjugacy classes of $Dic_n$ and the orders of their centralizer subgroups are given as follows:
\begin{equation}\label{conj_class}
\begin{split}
&\{1\}, \quad \{a^n\}\quad (\text{order of centralizer group}=4n)\\
&\{a^l, a^{-l}\}, 1\leq l\leq n-1,\quad (\text{order of centralizer group}=2n)\\
&\{xa,xa^3,xa^{5},...,xa^{2n-1}\}, \quad \{x, xa^2, xa^4,...,xa^{2n-2}\} \quad (\text{order of centralizer group}=4).
\end{split}
\end{equation}
Using (\ref{action_typeD}) and (\ref{conj_class}) it is easy to identify the contribution from each conjugacy class. It follows that (\ref{twisted_sector_contribution}) is given by 
\begin{equation}
\frac{1}{2n}\sum_{k=1}^{n-1}\frac{1}{2-\exp(\pi\sqrt{-1}k/n)-\exp(\pi\sqrt{-1}k/n)^{-1}}+\frac{1}{16n}+\frac{1}{8}+\frac{1}{8}.
\end{equation}
We need to evaluate the sum $\sum_{k=1}^{n-1}\frac{1}{2-\exp(\pi\sqrt{-1}k/n)-\exp(\pi\sqrt{-1}k/n)^{-1}}$. Again by \cite[Lemma 3.3.2.1]{lie}, we have 
\begin{equation*}
\begin{split}
\frac{(2n)^2-1}{12}
=&\sum_{k=1}^{2n-1}\frac{1}{2-\exp(2\pi\sqrt{-1}k/(2n))-\exp(2\pi\sqrt{-1}k/(2n))^{-1}}\\
=&\sum_{k=1}^{n-1}\frac{1}{2-\exp(\pi\sqrt{-1}k/n)-\exp(\pi\sqrt{-1}k/n)^{-1}}+\frac{1}{4}\\
&+ \sum_{k=1}^{n-1}\frac{1}{2-\exp(2\pi\sqrt{-1}(n+k)/(2n))-\exp(2\pi\sqrt{-1}(n+k)/(2n))^{-1}}.
\end{split}
\end{equation*}
Note that
\begin{equation*}
\begin{split}
&2-\exp(2\pi\sqrt{-1}(n+k)/(2n))-\exp(2\pi\sqrt{-1}(n+k)/(2n))^{-1}\\
=&2+\exp(\pi\sqrt{-1}k/n)+\exp(\pi\sqrt{-1}k/n)^{-1}\\
=&2+2\cos(\pi k/n)=4\cos^2(\pi k/(2n))=4\sin^2((\pi(k+n)/(2n));\\
&2-\exp(\pi\sqrt{-1}k/n)-\exp(\pi\sqrt{-1}k/n)^{-1}\\
=&2-2\cos(\pi k/n)=4\sin^2(\pi k/(2n)).
\end{split}
\end{equation*}
Since $\sin(\pi (k+n)/(2n))=-\sin(\pi(k-n)/(2n))$, we see that 
\begin{equation*}
\begin{split}
&\sum_{k=1}^{n-1}\frac{1}{2-\exp(\pi\sqrt{-1}k/n)-\exp(\pi\sqrt{-1}k/n)^{-1}}\\
&= \sum_{k=1}^{n-1}\frac{1}{2-\exp(2\pi\sqrt{-1}(n+k)/(2n))-\exp(2\pi\sqrt{-1}(n+k)/(2n))^{-1}},
\end{split}
\end{equation*}
 from which it follows that $$2 \sum_{k=1}^{n-1}\frac{1}{2-\exp(\pi\sqrt{-1}k/n)-\exp(\pi\sqrt{-1}k/n)^{-1}}+\frac{1}{4}=\frac{(2n)^2-1}{12}.$$ This shows that $$\sum_{k=1}^{n-1}\frac{1}{2-\exp(\pi\sqrt{-1}k/n)-\exp(\pi\sqrt{-1}k/n)^{-1}}=\frac{n^2-1}{6}$$ and (\ref{twisted_sector_contribution}) is given by 
 $$\frac{n^2-1}{12n}+\frac{1}{16n}+\frac{1}{8}+\frac{1}{8}=\frac{1}{12}(n+3-\frac{1}{4n}).$$
 Since the exceptional divisor of the minimal resolution of $\mathbb{C}^2/Dic_n$ is a tree of $\mathbb{CP}^1$ whose dual graph is the Dynkin diagram $D_{n+2}$, its Euler characteristic is $n+3$ and the Lemma is proved in this case.
 
\subsection{Type E}
If $\mathbb{C}^2/G_i$ is of type $E$, then there are three possibilities: $E_6, E_7, E_8$. The group $G_i$ is isomorphic to the binary tetrahedral group (for $E_6$), the binary octahedral group (for $E_7$), or the binary icosahedral group (for $E_8$). In each case the group and its action on $\mathbb{C}^2$ can be explicitly described, and the Lemma can be proved by computing (\ref{twisted_sector_contribution}) using this information. We work out the details for $E_6$ and leave the other two cases to the reader.

In the $E_6$ case, the group $G_i$ is isomorphic to the binary tetrahedral group $2T$. This group is of order $24$ and its elements can be identified with the following quaternion numbers: $$\frac{1}{2}(\pm 1\pm i\pm j\pm k), \quad \pm i, \quad \pm j, \quad \pm k, \quad, \pm 1.$$ The group $2T$ has $7$ conjugacy classes:\\\\
\begin{tabular}[b]{|c|c|c|c|c|}
\hline
Conjugacy Class & $(1)$ & $(-1)$ & $(i)$ & $(\frac{1}{2}(1+i+j+k))$\\ 
\hline
Size & $1$ & $1$ & $6$ & $4$\\ 
\hline
Conjugacy Class & $(\frac{1}{2}(1+i+j-k))$ & $(\frac{1}{2}(-1+i+j+k))$ & $(\frac{1}{2}(-1+i+j-k))$ & \\
\hline
Size & $4$ & $4$ & $4$ & \\
\hline
\end{tabular}
The action of $2T$ on $\mathbb{C}^2$ can be described using the following identification $$x+yi+zj+wk\mapsto \left(\begin{array}{cc} x+yi& z+wi\\ -z+wi& x-yi\end{array}\right).$$
Now it is straightforward to see that (\ref{twisted_sector_contribution}) is given by $$\frac{1}{24}\frac{1}{2-(-2)}+\frac{1}{4}\frac{1}{2-0}+\frac{1}{6}\frac{1}{2-1}+\frac{1}{6}\frac{1}{2-1}+\frac{1}{6}\frac{1}{2-(-1)}+\frac{1}{6}\frac{1}{2-(-1)}=\frac{167}{288}=\frac{1}{12}(7-\frac{1}{24}).$$
Since $7$ is the Euler characteristic of the exceptional divisor of the minimal resolution of $\mathbb{C}^2/2T$, the result follows.

\end{document}